\input amstex
\documentstyle{amsppt}
\magnification=\magstep1
\hsize=5.4in
\vsize=6.8in
\topmatter

\centerline  {\bf RIGIDITY FOR EQUIVALENCE RELATIONS ON HOMOGENEOUS SPACES}

\vskip .15in
\centerline {\rm ADRIAN IOANA\footnote{ Supported by a Clay Research Fellowship.} and YEHUDA SHALOM\footnote{Supported by NSF grants DMS 0701639 and 1007227.}}

\address Mathematics Department,  UCLA, Los Angeles, CA 91125
\endaddress
\email adiioana\@math.ucla.edu \endemail

\address Mathematics Department,  UCLA, Los Angeles, CA 91125
\endaddress
\email shalom\@math.ucla.edu\endemail

\abstract 
We study Popa's notion of rigidity for equivalence relations induced by actions on homogeneous spaces.
 For any lattices $\Gamma,\Lambda$ in a semisimple Lie group $G$ with finite center and no compact factors we prove that  
  the   action $\Gamma\curvearrowright G/\Lambda$ is rigid. 
If in addition $G$ has property (T) then we derive that the von Neumann algebra $L^{\infty}(G/\Lambda)\rtimes\Gamma$ has property (T).
We also show that if the adjoint action of $G$ on the Lie algebra of $G$ $-$ $\{0\}$ is amenable (e.g. if $G=SL_2(\Bbb R)$), then
 any ergodic subequivalence relation of the orbit equivalence relation 
 of the action  $\Gamma\curvearrowright G/\Lambda$ is either hyperfinite or rigid.

\endabstract 
\endtopmatter

\document

\head Introduction and statement of main results.\endhead
\vskip 0.1in

In [Po06], S. Popa introduced the notion of {\it relative property (T)}   
for inclusions of von Neumann algebras $B\subset M$ -- it has since been at the heart of his deformation/rigidity theory. 
When applied to inclusions arising from actions and equivalence relations, this concept suggested two new properties for actions
and equivalence relations:
\vskip 0.03in

\hskip 0.05in$\bullet$ A probability measure preserving (pmp) action $\Gamma\curvearrowright (X,\mu)$ of a countable group  $\Gamma$  

\hskip 0.17in is {\it rigid} if the inclusion of $L^{\infty}(X)$ in the crossed-product algebra   $L^{\infty}(X)\rtimes\Gamma$ ([MvN36])

\hskip 0.17in   has the relative property (T) in the sense of [Po06, Definition 4.2.1].

\hskip 0.05in $\bullet$ A countable pmp equivalence relation $R$ on $(X,\mu)$ is {\it rigid} if the inclusion of $L^{\infty}(X)$ 

\hskip 0.17in in the von Neumann algebra $L(R)$ of $R$ ([FM77]) has the relative property (T).
\vskip 0.03in

Note that for free actions, rigidity is a property of their equivalence relations: a {\it free} pmp action $\Gamma\curvearrowright (X,\mu)$ is rigid if and only if its {\it orbit equivalence  relation} ($x\sim y$ if $\Gamma x=\Gamma y$) is.

\vskip 0.05in
In the last decade these notions of rigidity have led to several remarkable applications, most notably, to calculations of invariants of von Neumann algebras and equivalence relations ([Po06],[PV10],[Ga08]) and to constructions of non-orbit equivalent actions of non-amenable groups ([GP05],[Io07]).

Yet, while rigidity was successfully exploited in applications, the theoretical aspects of its study (e.g. finding new constructions of rigid equivalence relations and a more manageable definition of rigidity -- avoiding the use of von Neumann algebras)
were neglected. In fact, until recently all known examples of rigid actions and equivalence relations ([Po06],[Ga08]) relied on the following group theoretic construction. Let $A$ be a countable abelian group together with an action of a countable group $\Gamma$ such that the pair $(\Gamma\ltimes A,A)$ has the relative property (T) of Kazhdan-Margulis. Then the (Haar) measure preserving
 action $\Gamma\curvearrowright \hat{A}$ and its orbit equivalence relation are rigid ([Po06]). In particular, since the pair (SL$_n(\Bbb Z)\ltimes\Bbb Z^n,\Bbb Z^n)$ has the relative property (T) ([Ka67],[Ma82]), it follows that the natural action of SL$_n(\Bbb Z)$ on the $n$-torus $\Bbb T^n$ is rigid for all $n\geqslant 2$.

The situation improved with the finding of an ergodic theoretic criterion for rigidity of pmp actions and equivalence relations ([Io10], see the end of the introduction). The criterion was then used to produce  the first examples of rigid equivalence relations not built from a pair of groups with relative property (T): if $S$ denotes the orbit equivalence relation of the action SL$_2(\Bbb Z)\curvearrowright \Bbb T^2$, then any ergodic non-hyperfinite subequivalence relation $R\subset S$ is  rigid ([Io10, Theorem 0.1]). 
Although this result provides new instances of rigidity, it has the disadvantage of being limited to a specific action. 

In this paper, we work in the general framework of actions on homogeneous spaces and prove rigidity for the induced (sub)equivalence relations under fairly general assumptions. More precisely, we consider actions of countable subgroups $\Gamma<G$ on the homogeneous space $(G/\Lambda,m_{G/\Lambda})$, where $G$ is a real algebraic group, $\Lambda<G$ is a lattice and $m_{G/\Lambda}$ is the unique $G$-invariant probability measure on $G/\Lambda$.

Our first result asserts that, under mild assumptions on $G$, the action of any lattice $\Gamma<G$ on $G/\Lambda$ is rigid.

\proclaim {Theorem  A} Let $G$ be a real algebraic group with finite center, no proper normal co-compact algebraic subgroups, and no non-trivial algebraic homomorphism into $\Bbb R^*$.  
\vskip 0.05in
\noindent
If $\Gamma,\Lambda<G$ are lattices, then the pmp action $\Gamma\curvearrowright (G/\Lambda,m_{G/\Lambda})$ is rigid.

\noindent
Moreover, if $G$ has property (T) (e.g. if $G$ is a connected semisimple Lie group with finite center whose simple factors have real-rank $\geqslant$ 2), then  $L^{\infty}(G/\Lambda)\rtimes\Gamma$ has property (T).
\endproclaim

Property (T) for von Neumann algebras $M$ was introduced by Connes and Jones in [CJ85]. For a crossed-product algebra $L^{\infty}(X)\rtimes\Gamma$ coming from a pmp action, it is equivalent to having both that $\Gamma$ is a property (T) group and that the action $\Gamma\curvearrowright (X,\mu)$ is rigid. Because lattices inherit property (T) ([Ka67]), this indicates how to deduce the last line of Theorem A.

Theorem A implies that for any $n\geqslant 3$, the crossed product von Neumann algebra $L^{\infty}($SL$_n(\Bbb R)/$SL$_n(\Bbb Z))\rtimes $SL$_n(\Bbb Z)$ has property (T). It is worth mentioning that these  are the first examples of property (T) von Neumann algebras that are not constructed from countable property (T) groups.

As a consequence of Theorem A we also derive:

\proclaim {Corollary B}
Let $G$ be as in Theorem A and $\Gamma,\Lambda<G$ be lattices. Let $I$ be a countable set on which $\Gamma$ acts with  infinite orbits and $X_0$ be some probability space. Endow $X=X_0^I$  with the corresponding generalized Bernoulli $\Gamma$-action.

\vskip 0.05in
\noindent
 Then for any  pmp $\Gamma$-space $Y$, any measurable quotient $\Gamma$-map $p:X\times Y\rightarrow G/\Lambda$ depends a.e. on the second coordinate only. 
\endproclaim

The main result of [Io10] shows that the orbit equivalence relation $S$ of the action SL$_2(\Bbb Z)\curvearrowright \Bbb T^2$ satisfies the following ``global" dichotomy: any ergodic subequivalence relation $R\subset S$ is either hyperfinite or rigid. 
Our second result establishes this dichotomy for many other actions, including the action of SL$_2(\Bbb Z)$ on SL$_2(\Bbb R)/$SL$_2(\Bbb Z)$:

\proclaim {Theorem C} Let $G$ be any of the groups SL$_2(\Bbb R)$, SL$_2(\Bbb C)$, SL$_2(\Bbb R)\ltimes\Bbb R^2$ or SL$_2(\Bbb C)\ltimes\Bbb C^2$. Let $\Gamma<G$ be a countable discrete subgroup and $\Lambda<G$ be a lattice. Denote by $S$ the orbit equivalence relation of the action $\Gamma\curvearrowright (G/\Lambda,m_{G/\Lambda})$.
\vskip 0.05in
\noindent
Then any ergodic subequivalence relation $R\subset S$ is either hyperfinite or rigid. 

\noindent
Moreover, for any subequivalence relation $R\subset S$, we can find a measurable partition $G/\Lambda=X_0\cup X_1$ such that $X_0,X_1$ are $R$-invariant, $R_{|X_0}$ is hyperfinite and $R_{|X_1}$ is rigid.

\endproclaim
\noindent
The rest of the paper consists of two sections. In the next one we recall two ergodic theoretic criteria 
for rigidity of actions and equivalence relations.
In the last section, we use these criteria to prove Theorems A and C, and their more general versions, Theorems D and E.
\vskip 0.1in
\noindent
{\it Acknowledgment.} We would like to express our deep gratitude to Gregory Margulis who made a crucial contribution to the paper during his distinguished lecture series visit at UCLA. 
\vskip 0.2in
\head criteria for rigidity\endhead

\vskip 0.2in

Theorem 4.4 in [Io10] gives an ergodic theoretic formulation of rigidity for {\it free ergodic} actions. Also, Proposition 2.2 in [Io10] provides an ergodic theoretic criterion for rigidity of {\it ergodic} equivalence relations. In this section, we note that appropriate versions of these results
 -- implicitly proved, but not stated in [Io10] --  hold without the freeness and ergodicity assumptions.

\proclaim {\bf Proposition 1 (equivalent formulation of rigidity for actions)}

\noindent
A pmp action $\Gamma\curvearrowright (X,\mu)$ of a countable group $\Gamma$ on a probability space $(X,\mu)$
is rigid if and only if for any sequence of Borel probability measures $\nu_n$ on $X\times X$ satisfying:

\hskip 0.05in {\bf (1)} $p_*^i\nu_n=\mu$ for all $n$ and $i=1,2$, where $p^i:X\times X\rightarrow X$ denotes the projection onto 

\hskip 0.3in the $i$-th coordinate.

\hskip 0.05in {\bf (2)} $\int_{X\times X}f(x)g(y)d\nu_n(x,y)\rightarrow \int_{X}f(x)g(x)d\mu(x)$, for all bounded Borel functions 

\hskip 0.3in $f,g:X\rightarrow \Bbb C$.

\hskip 0.05in {\bf (3)} $||(\gamma\times\gamma)_*\nu_n-\nu_n||\rightarrow 0$, for every $\gamma\in\Gamma$.

\noindent
we have that $\nu_n(\Delta)\rightarrow 1$ (where $\Delta\subset X\times X$ denotes the diagonal).
\endproclaim

\noindent
Here, for a bounded signed Borel measure $\nu$ on $X\times X$, the norm $||\nu||$ is obtained by viewing $\nu$ as a linear functional on the space of bounded Borel functions on $X\times X$.

\proclaim 
{\bf Proposition 2 (criterion for rigidity of equivalence relations)} 

\noindent
Let $R$ be a countable pmp equivalence relation on a probability space $(X,\mu)$. 
Assume that for any sequence of Borel probability measures $\nu_n$ on $X\times X$ satisfying {\bf (1)}, {\bf (2)} and

\hskip 0.05in {\bf (3')} $||(\theta\times\theta)_*\nu_n-\nu_n||\rightarrow 0$, for every $\theta$ belonging to the group $[R]$ of automorphisms

\hskip 0.35in  of $(X,\mu)$ whose graph is contained in $R$.

\noindent
we have that $\nu_n(\Delta)\rightarrow 1$.

\vskip 0.05in
\noindent Then $R$ is rigid.

\endproclaim

\vskip 0.05in

Before indicating how these propositions follow from [Io10], let us recall the notion of relative property (T) for von Neumann algebras.

\vskip 0.05in
\noindent
{\bf Definition [Po06, Definition 4.2.1]}. Let $(M,\tau)$ be a von Neumann algebra with a normal faithful tracial state $\tau$ and $B\subset M$ a von Neumann subalgebra.

\noindent
We say that  $B\subset M$ has {\it relative property (T)} if whenever $\Cal H$ is a Hilbert $M$-bimodule and $\xi_n\in\Cal H$ is a sequence satisfying 

\hskip 0.05in $\bullet$ $\langle x\xi_n,\xi_n\rangle=\langle\xi_n x,\xi_n\rangle=\tau(x)$, for all $x\in M$ and every $n\geqslant 1$. \hskip 0.15in ({\it tracial})

\hskip 0.05in $\bullet$ $||x\xi_n-\xi_n||\rightarrow 0$, for all $x\in M$. \hskip 0.15in({\it almost central})

\noindent
we can find $\eta_n\in \Cal H$ such that $b\eta_n=\eta_nb$, for all $b\in B$ and every $n\geqslant 1$, and
$||\eta_n-\xi_n||\rightarrow 0$.

\vskip 0.05in

\noindent {\it Proof of Proposition 1.} The proof of the ``only if part" is identical to that of ``(a) $\Rightarrow$ (c)"  in [Io10, Theorem 4.4] (which does not actually use the freeness and ergodicity assumptions).

Now, let $\Gamma\curvearrowright (X,\mu)$ be an  action such that for any sequence of probability measures $\nu_n$ satisfying {\bf (1)}-{\bf (3)}, we must have $\nu_n(\Delta)\rightarrow 1$.	Denote $M=L^{\infty}(X)\rtimes\Gamma$. To prove that the action is rigid, let $\Cal H$ be a Hilbert $M$-bimodule and $\xi_n\in\Cal H$ a sequence of tracial, almost central vectors. Denote by $\Cal H_n\subset\Cal H$ the cyclic $L^{\infty}(X)$-bimodule $\overline{L^{\infty}(X)\xi_nL^{\infty}(X)}$.

 Let $\nu_n$ be a probability measure on $X\times X$ such that $\Cal H_n\ni f\xi_n g\rightarrow f\otimes g\in L^2(\nu_n)$ extends to an isomorphism of $L^{\infty}(X)$-bimodules (where $(f\otimes g)(x,y)=f(x)g(y)$).
Then the proof Lemma 2.1 in [Io10] implies that $\nu_n$ satisfy conditions {\bf (1)}-{\bf (3)}. Thus, we must have $\nu_n(\Delta)\rightarrow 1$.
Since $\eta_n=1_{\Delta}\in L^2(\nu_n)\cong \Cal H_n$ verifies $f\eta_n=\eta_nf$, for all $f\in L^{\infty}(X)$, and $||\eta_n-\xi_n||_2=\sqrt{1-\nu_n(\Delta)}\rightarrow 0,$ we are done. \hfill$\square$
\vskip 0.05in
\noindent
{\it Proof of Proposition 2.} This is the same as the proof of the ``if part" of Proposition 1, where we replace $M$ by $L(R)$.
\hfill$\square$

\vskip 0.2in
\head Proofs\endhead

\vskip 0.2in
We begin by stating the more general versions of Theorems A and C. 
\proclaim {Theorem D} Let $G$ be a real algebraic group and  $\Lambda\subset G$ be a lattice. 
Let $\Gamma\subset G$ be a  countable subgroup and denote by $H$ its Zariski closure. Assume that $H$ has no proper normal co-compact algebraic subgroup and no non-trivial homomorphism into $\Bbb R^*$.  
Let $\eta$ be a $\Gamma$-invariant probability measure on $G/\Lambda$. 

\vskip 0.05in
\noindent  
If the centralizer of $\Gamma$ (equivalently, of $H$) in $G$ is finite, then the pmp action $\Gamma\curvearrowright (G/\Lambda,\eta)$ is rigid. 

\vskip 0.01in
\noindent
In the case $\eta=m_{G/\Lambda}$,  the converse is true: if the action $\Gamma\curvearrowright (G/\Lambda,\eta)$ is rigid, then the centralizer of $\Gamma$ in $G$ is finite.
\endproclaim 

{\it Remark.} Theorem D implies that for $\Gamma=\Bbb F_2\times\Bbb Z$ actions of the form $\Gamma\curvearrowright G/\Lambda$ are never rigid. It would be interesting to decide whether $\Gamma$ admits a free ergodic rigid action at all. Note in this respect that the general question of characterizing non-amenable groups which admit free ergodic rigid actions ([Po06, Problem 5.10.2]) remains open (see [Ga08] for a partial result).

\proclaim {Theorem E} Let $G$  be a real algebraic group, $m=m_{G}$ be a the Haar measure on $G$ and on $(G,m)$ consider the left-right multiplication action of $G\times G: (g_1,g_2)\cdot g=g_1gg_2^{-1}$. 

\noindent
Let $\Gamma,\Lambda<G$ be two countable discrete subgroups and denote by $H, K$ their Zariski closures.
Assume that the stabilizer  of any point in the Lie algebra of $G$ under the adjoint actions of $H$ and $K$ is amenable.

\noindent
Denote by $S$ the orbit equivalence relation of the action $\Gamma\times\Lambda\curvearrowright (G,m)$, i.e. $S=\{(x,y)\in G\times G|x\in \Gamma y\Lambda\}$. Let $X\subset G$ be a Borel set with $0<m(X)<\infty$.

\vskip 0.05in
\noindent
Then for any subequivalence relation $R\subset S_{|X}=S\cap (X\times X)$, we can find an $R$-invariant measurable partition $X=X_0\cup X_1$ such that $R_{|X_0}$ is hyperfinite and $R_{|X_1}$ is rigid.

\endproclaim

{\it Remark.} The assumption that $\Gamma,\Lambda<G$ are discrete is essential. Otherwise, we would allow the case $G$ is compact. If $G$ is compact, then we can find a $G\times G$-invariant metric $d$ on $G$ (defining the topology).  For $n\geqslant 1$, let $A_n=\{(x,y)\in G\times G|d(x,y)<\frac{1}{n}\}$ and set $\nu_n=\frac{(m\times m)_{|A_n}}{(m\times m)(A_n)} $. Then $\nu_n$ is a  Borel probability measure on $G\times G$ which is invariant under the diagonal product action of $G\times G$ on $G\times G$: $(g_1,g_2)\cdot (h,k)=(g_1hg_2^{-1},g_1kg_2^{-1})$ and whose projection onto both coordinates is equal to $m$. Moreover, $\nu_n$  converge weakly to the pushforward of $m$ 
through the map $G\ni x\rightarrow (x,x)\in G\times G$. This shows that the  action  $\Gamma\times\Lambda\curvearrowright (G,m)$ is not rigid, for countable subgroups $\Gamma,\Lambda<G$.
\vskip 0.1in 
\noindent
Next, we introduce some notations that we will use in the proofs of Theorems D and E. Fix an unimodular real algebraic group  $G$ and a Haar measure $m=m_{G}$. 

\hskip 0.05in $\bullet$ Denote by g the Lie algebra of $G$ and by P(g)=(g$\setminus\{0\})/\Bbb R^{*}$ the associated projective

\hskip 0.17in variety together with the map $p:\text{g}\setminus\{0\}\rightarrow$ P(g).

\hskip 0.17in  We endow  g and P(g) with the 
adjoint $G$-action: Ad$(\gamma)(x)=\gamma x{\gamma}^{-1}$.

\hskip 0.05in $\bullet$ Let $q: G\rightarrow$ g be a Borel map which is equal to the matrix logarithm in some 

\hskip 0.17in  neighborhood  $U$ of $1\in G$.

\hskip 0.05in $\bullet$ Next, define $r:G\times G\rightarrow G$ by $r(x,y)=xy^{-1}$.

\hskip 0.05in $\bullet$ We can now set $\rho=p\circ q\circ r:(G\times G)\setminus \Delta\rightarrow$ P(g), where $\Delta=\{(x,x)|x\in G\}$.

\hskip 0.05in $\bullet$ Finally, we let $\pi:(G\times G)\setminus\Delta\rightarrow G\times$P(g) be given by $\pi(x,y)=(x,\rho(x,y))$.
\vskip 0.03in
\noindent
Given a Borel subset $X\subset G$ we denote $\Delta_X=\Delta\cap (X\times X)$ and by $p^i:X\times X\rightarrow X$ the projection onto the $i$-th coordinate.
 For a Polish space $Y$, we denote by $\Cal M(Y)$ the space of Borel probability measures on $Y$.

\proclaim {Lemma F} Let $X\subset G$ be a Borel subset endowed with a Borel probability measure $\eta$. 

\noindent
Let $c>0$. Suppose that $\mu_n\in\Cal M(X\times X)$ is a sequence  satisfying  $p^1_{*}\mu_n\leqslant c\eta$, $\mu_n(\Delta_X)=0$, for all  $n\geqslant 1$, and $\mu_n(A\times (X\setminus A))\rightarrow 0$, for any Borel set $A\subset X$.

\noindent
Let $\Gamma\subset G$ be a countable subgroup and $\phi_1,\phi_2:X\rightarrow \Gamma$ be two Borel maps.  Denote by $D$ the set of $(x,y)\in (X\times X)\setminus \Delta_X$ such that $\rho(\phi_1(x)x\phi_2(x),\phi_1(y)y\phi_2(y))=\text{Ad}(\phi_1(x))(\rho(x,y))$.

\vskip 0.05in
\noindent
Then $\lim_{n\rightarrow\infty}\mu_n(D)=1.$ 

\endproclaim

\noindent
{\it Proof.} We first claim that if $\{B_i\}_{i=1}^{\infty}$ is a Borel partition of $X$, then $\mu_n(\cup_{i=1}^{\infty}(B_i\times B_i))\rightarrow 1$. For $k\geqslant 1$, set $X_k=\cup_{i=1}^kB_i$. Note that $(X\times X)\setminus\cup_{i=1}^{\infty}(B_i\times B_i)\subset \cup_{i=1}^k(B_i\times (X\setminus B_i))\cup ((X\setminus X_k)\times X)$. Since $\mu_n(B_i\times (X\setminus B_i))\rightarrow 0$, for all $i$,  and $p^1_{*}\mu_n\leqslant c\eta$, we deduce that $$\limsup \mu_n((X\times X)\setminus\cup_{i=1}^{\infty}(B_i\times B_i))\leqslant c\eta(X\setminus X_k),\forall k\geqslant 1.$$

Since $\{B_i\}_{i=1}^{\infty}$ is a partition of $X$, we get that $\eta(X_k)\rightarrow 1$, which proves our claim.

Towards proving that  $\mu_n(D)\rightarrow 1$, let $\varepsilon>0$. Then we can find a finite subset $F$ of $\Gamma$ such that $\eta(\{x\in X|\phi_1(x)\in F\})\geqslant 1-\frac{\varepsilon}{c}.$
Since $p^1_*\mu_n\leqslant c\eta$ it follows that 
 $A=\{(x,y)\in X\times X|\phi_1(x)\in F\}$ satisfies $\mu_n(A)\geqslant 1-\varepsilon,$ for all $n\geqslant 1$.

Next, since $q(x)=\log(x)$, for $x$ in a neighborhood of $1\in G$, we get that for all $x$ in a, possibly smaller, neighborhood $V$ of $1$ we have that $q(\gamma x\gamma^{-1})=\gamma q(x)\gamma^{-1}$, for every $\gamma\in F$.
Let $B=\{(x,y)\in X\times X|xy^{-1}\in V\}$. Let $W$ be a neighborhood of $1\in G$ such that $WW^{-1}\subset V$ and $h_1,h_2,..\in G$ be a sequence such that $G=\cup_{i=1}^{\infty}Wh_i$. For $i\geqslant 1$, define $B_i=(Wh_i\setminus(\cup_{j=1}^{i-1}Wh_j))\cap X$. Since  $\{B_i\}_{i=1}^{\infty}$ is a partition of $X$  and $\cup_{i=1}^{\infty}(B_i\times B_i)\subset B$, the above claim yields that $\mu_n(B)\rightarrow 1$.

Finally, let us show that $C=\{(x,y)\in X\times X|\phi_1(x)=\phi_1(y),\phi_2(x)=\phi_2(y)\}$ satisfies $\mu_n(C)\rightarrow 1$. This also follows for the above claim, after noticing that the sets $C_{\gamma_1,\gamma_2}=\{x\in X|\phi_1(x)=\gamma_1,\phi_2(x)=\gamma_2\}$, with $\gamma_1,\gamma_2\in\Gamma$, form a partition of $X$
and satisfy $C=\cup_{\gamma_1,\gamma_2\in\Gamma}(C_{\gamma_1,\gamma_2}\times C_{\gamma_1,\gamma_2})$. Finally, it is easy to see that if $A\cap B\cap C\subset D$. Since $\liminf\mu_n(A\cap B\cap C)\geqslant 1-\varepsilon$ and $\varepsilon>0$ is arbitrary, the proof is complete.\hfill$\square$

\vskip 0.1in

\noindent
{\it Proof of Theorem D.} Suppose first that the centralizer of $H$ in $G$ is finite. Let $X\subset G$ be a fundamental domain for the right $\Lambda$-action. Identify $G/\Lambda$ with $X$ via the map $G/\Lambda\ni x\Lambda\rightarrow x\Lambda\cap X\in X$. Under this identification, the corresponding $\Gamma$-action on $X$ is  given by $\gamma\cdot x=\gamma x w(\gamma,x)$, for all $x\in X$ and $\gamma\in\Gamma$, where $\lambda=w(\gamma,x)$ is the unique element of $\Lambda$ such that $\gamma x\lambda\in X$.

Let $\nu_n$ be a sequence of Borel probability measures on $X\times X$ satisfying
\vskip 0.02in
\noindent
 (1) $p_*^i\nu_n=\eta$, for all $n$ and $i=1,2$. 

\noindent (2) $\int_{X\times X}f(x)g(y)d\nu_n(x,y)\rightarrow\int_{X}fgd\eta$, for all bounded Borel functions $f,g:X\rightarrow\Bbb C$. 

\noindent (3) $||(\gamma\times\gamma)_*\nu_n-\nu_n||\rightarrow 0$, for all $\gamma\in\Gamma$. 

\vskip 0.05in

By Proposition 1, to conclude that the action $\Gamma\curvearrowright (X,\eta)$ is rigid, it suffices to argue that $\nu_n(\Delta_X)\rightarrow 1$. 

If this is false, then after passing to a subsequence we may assume that $c_n=1-\nu_n(\Delta_X)$ verify $c=\inf c_n>0$. Define $\mu_n\in\Cal M(X\times X)$ by $\mu_n(A)=c_n^{-1}\nu_n(A\setminus\Delta_X)$, for any Borel set $A\subset X\times X$. Then conditions (1) and (2) imply that $p_*^1\mu_n\leqslant c^{-1}\eta$, for all $n$, and  $\mu_n(A\times (X\setminus A))\rightarrow 0$, for any Borel set $A\subset X$. By applying Lemma F, we get that $\mu_n(\{(x,y)\in (X\times X)\setminus\Delta_X|\rho(\gamma\cdot x,\gamma\cdot y)=\text{Ad}(\gamma)(\rho(x,y))\})\rightarrow 1$.

Also, condition (3) gives that $||(\gamma\times\gamma)_*\mu_n-\mu_n||\rightarrow 0$, for all $\gamma\in\Gamma$.
Combining the last two facts yields that the probability measures $\zeta_n=\rho_*\mu_n$ on $P$(g) satisfy $||\text{Ad}(\gamma)_*\zeta_n-\zeta_n||\rightarrow 0$, for all $\gamma\in\Gamma$. By taking a weak limit, we find a probability measure $\zeta$ on $P$(g) which is invariant under the adjoint action of $\Gamma$. 

By applying [Sh99, Theorem 3.11] we  get that $\zeta$ is invariant under the adjoint action of $H$ and that  $H$ has a normal co-compact subgroup which fixes every point in the support of $\zeta$. The hypothesis forces that $H$ fixes every point in the support of $\zeta$.
  In particular, there exists  $x\in$ g$\setminus\{0\}$ and a homomorphism $\chi:H\rightarrow\Bbb R^*$ such that $\gamma x\gamma^{-1}=\chi(\gamma) x$, for all $\gamma\in H$. Since every such $\chi$ is trivial, we deduce that $x$ commutes with $H$. Hence, for all $n$, we have that $x_n=$ exp$(\frac{x}{n})\in G$ commutes with $H$.  Since $x\not=0$, this contradicts the assumption that the centralizer of $H$ in $G$ is finite.

For the converse, suppose that $\eta=m_{G/\Lambda}$ and that the action $\Gamma\curvearrowright (G/\Lambda,\eta)$ is rigid.  By contradiction, if the centralizer of $H$ in $G$ is infinite, we can find a sequence of elements $x_n\in G\setminus\{1\}$ which commute with $H$ and converge to 1. For every $n$, let  $\nu_n$ be the pushforward of $\eta$ through the map $G/\Lambda\ni g\Lambda\rightarrow (g\Lambda,x_ng\Lambda)\in G/\Lambda\times G/\Lambda$. Since $x_n$ and $\Gamma$ commute, $\nu_n$ is invariant under the diagonal $\Gamma$-action. It is clear that $\nu_n$ converge weakly to the pushforward of $\eta$ through the map $G/\Lambda\ni x\rightarrow (x,x)\in G/\Lambda\times G/\Lambda$ and that the projection of $\nu_n$ onto both coordinates is equal to $\eta$. Since the action $\Gamma\curvearrowright (G/\Lambda,\eta)$ is  rigid, we conclude that $\nu_n(\Delta)\rightarrow 1$. Equivalently, $m(\{(g\in X|g^{-1}x_ng\notin\Lambda\})\rightarrow 0$, for every Borel set $X\subset G$ with $m(X)<\infty$. On the other hand, as $x_n\rightarrow 1$,  we have that $m(\{(g\in X|g^{-1}x_ng\notin U\})\rightarrow 0$, for any neighborhood $U$ of $1$ in $G$. Since $\Lambda$ is discrete, we  deduce that $x_n=1$, for $n$ large enough, a contradiction.
\hfill $\square$
\vskip 0.1in

\noindent
{\it Proof of Theorem E.}  In order to get the conclusion it suffices to show that if $R\subset S_{|X}$ is a non-rigid subequivalence relation, then there exists an $R$-invariant Borel subset $X_0\subset X$ such that $m(X_0)>0$ and $R_{|X_0}$ is hyperfinite.
Let $\eta=m(X)^{-1}m_{|X}$ be the probability measure on $X$ obtained by normalizing the restriction of $m$ to $X$. Since $R$ is not rigid, Proposition 2 gives a sequence  $\nu_n$ of Borel probability measures on $X\times X$ such that 
\vskip 0.02in
\noindent (1) $p_*^i\nu_n=\eta$ for all $n$ and $i=1,2$. 

\noindent (2) $\int_{X\times X}f(x)g(y)d\nu_n(x,y)\rightarrow \int_{X}fgd\eta$, for all bounded Borel functions $f,g:X\rightarrow\Bbb C$. 

\noindent (3) $||(\theta\times\theta)_*\nu_n-\nu_n||\rightarrow 0$, for every $\theta\in [R]$.

\noindent (4) $\nu_n(\Delta_X)\not\rightarrow 1$.

\vskip 0.05in
After passing to a subsequence we can assume that  $c_n=1-\nu_n(\Delta_X)$ satisfy $c=\inf c_n>0$. Define $\mu_n\in\Cal M(X\times X)$ by $\mu_n(A)=c_n^{-1}\nu_n(A\setminus\Delta_X)$, for every Borel set $A\subset X\times X$. We have that

\vskip 0.02in
\noindent (1') $p_*^i\mu_n\leqslant c^{-1}\eta$ for all $n$ and $i=1,2$. 

\noindent (2') $\mu_n(A\times (X\setminus A))\rightarrow 0$, for any Borel set $A\subset X$.

\noindent (3') $||(\theta\times\theta)_*\mu_n-\mu_n||\rightarrow 0$, for every $\theta\in [R]$.

\noindent (4') $\mu_n(\Delta_X)=0$, for all $n$.

\vskip 0.05in

Let $\theta\in [R]$. After modifying $\theta$ on a null set, we can assume that $\theta(x)\in\Gamma x\Lambda$, for all $x\in X$. This allows us to define $w_{\theta}=(w_{\theta}^1,w_{\theta}^2):X\rightarrow \Gamma\times\Lambda$   through the formula  $\theta(x)=w_{\theta}(x)\cdot x=w_{\theta}^1(x)xw_{\theta}^2(x)^{-1}$, for  every $x\in X$. 
Let  $\hat{\theta}$ be the Borel isomorphism of $X\times$P(g) given by  $\hat{\theta}(x,y)=(\theta(x),$ Ad$(w_{\theta}^1(x))(y))$. For every $n$, let $\zeta_n=\pi_{*}\mu_n$. Let us prove that $$||\hat{\theta}_{*}\zeta_n-\zeta_n||\rightarrow 0\tag a$$

By Lemma F (which applies as (1'),(2') and (4') hold true), the set $D$ of $(x,y)\in X\times X$ such that $\rho(\theta(x),\theta(y))=\text{Ad}(w_{\theta}^1(x))(\rho(x,y))$ satisfies $\mu_n(D)\rightarrow 1$. Since $(x,y)\in D$ if and only if $(\hat{\theta}\circ\pi)(x,y)=(\pi\circ(\theta\times\theta))(x,y)$, condition (3') gives (a).

Now, since $\zeta_n$ are Borel probability measures on the locally compact metrizable space 
$X\times$P(g), we can find a subsequence $\{\zeta_{n_k}\}_{k\geqslant 1}$ and a positive Borel measure
$\zeta$ on $X\times$P(g)  such that $\int fd\zeta_{n_k}\rightarrow \int fd\zeta$, for every $f\in C_0(X\times$P(g)).
We claim that $$\hat{\theta}_*\zeta=\zeta\hskip 0.07in\text{and}\hskip 0.07in\zeta(X\times P(\text{g}))=1\tag b$$

Note first that $0\leqslant \zeta(X\times$P(g)$)\leqslant 1$. Let $\varepsilon>0$. Let $X_0\subset X$ be a closed subset such that $\eta(X\setminus X_0)\leqslant\varepsilon$ and ${w_{\theta}}_{|X_0}:X_0\rightarrow\Gamma\times\Lambda$ is continuous. Condition (1') implies that $\zeta_n((X\setminus X_0)\times$P(g$))=\mu_n((X\setminus X_0)\times X)\leqslant c^{-1}\eta(X\setminus X_0)\leqslant c^{-1}\varepsilon$, for every $n$. Thus $\zeta_n(X_0\times $P(g))$\geqslant 1-c^{-1}\varepsilon$ and since $X_0\times$P(g) is closed, we deduce that $\zeta(X_0\times$P(g))$\geqslant 1-c^{-1}\varepsilon$. As $\varepsilon>0$ is arbitrary, this proves the second assertion from (b).

By (a), in order to prove that $\hat{\theta}_*\zeta=\zeta$ it suffices to show that $\int (f\circ\hat{\theta})d\zeta_{n_k}\rightarrow \int (f\circ\hat{\theta})d\zeta$, for every $f\in C_0(X\times$P(g)) with $||f||_{\infty}\leqslant 1$. Note that the restriction of $\hat{\theta}$ to $X_0\times$P(g) is continuous. Let $h\in C_0(X\times$P(g)) such that $||h||_{\infty}\leqslant 1$ and $h_{|X_0\times\text{P(g)}}=(f\circ\hat{\theta})_{|X_0\times\text{P(g)}}$. Since $\zeta((X\setminus X_0)\times$P(g))$\leqslant c^{-1}\varepsilon$ and $\zeta_n((X\setminus X_0)\times$P(g))$\leqslant c^{-1}\varepsilon$, we have that $$|\int (f\circ\hat{\theta})d\zeta_{n_k}-\int (f\circ\hat{\theta})d\zeta|\leqslant 4c^{-1}\varepsilon+ |\int hd\zeta_{n_k}-\int hd\zeta|.$$

Since $\int hd\zeta_{n_k}\rightarrow\int hd\zeta$ and $\varepsilon>0$ is arbitrary, this concludes the proof of (b). 
\vskip 0.05in

Next, since $p_*^1\mu_n\leqslant c^{-1}\eta$ and $\pi(x,y)=(x,...)$, we deduce that the push forward of $\zeta_n$ onto the $X$-coordinate is $\leqslant c^{-1}\eta$. Thus, the push forward of $\zeta$ onto the $X$-coordinate, denoted $\tilde\eta$, satisfies $\tilde\eta\leqslant c^{-1}\eta$. By using (b), we get that $\tilde\eta$ is $R$-invariant. 
Disintegrate $\zeta=\int_{X}\zeta_x d\tilde\eta(x)$, where $\zeta_x\in\Cal M($P(g)), for all $x\in X$. By using (b) and the fact that $\tilde\eta$ is $R$-invariant, the uniqueness of the disintegration implies that $$\zeta_{\theta(x)}={\text{Ad}(w_{\theta}^1(x))}_{*}\zeta_x, \hskip 0.07in \text{for}\hskip 0.05in \tilde\eta-\text {a.e.}\hskip 0.05in x\in X,\hskip 0.07in\forall \theta\in [R].$$

Denote by $X_0$ the support of $X$ and notice that it is $R$-invariant. Since $\tilde\eta\leqslant c^{-1}\eta$ we get that $\eta(X_0)>0$ and 

$$\zeta_{\theta(x)}={\text{Ad}(w_{\theta}^1(x))}_{*}\zeta_x, \hskip 0.07in \text{for}\hskip 0.05in \eta-\text {a.e.}\hskip 0.05in x\in X_0,\hskip 0.07in\forall \theta\in [R]\tag c$$

In the second half of the proof, we use (c) (and an analogous identity for $w_{\theta}^2$) to deduce that $R_{|X_0}$ is hyperfinite. We first do this under the additional assumption that $R_{|X_0}$ is {\it ergodic} with respect to $\eta$.

Since the adjoint action of $H$ on g is linear, by [Zi84, Corollary 3.2.12] the action of $H$ on $\Cal M$(P(g)) is smooth. Since $R_{|X_0}$ is ergodic, by using (c), we deduce that $\zeta_x$ lies in a single $H$-orbit, on a co-null subset of $X_0$. In other words, there exists $\xi\in\Cal M$(P(g)) such that $\zeta_x\in H\xi$, for $\eta$-almost every $x\in X_0$. Identify $H\xi$ with $H/P$, where $P$ denotes the stabilizer of $\xi$ in $H$.

We claim that $P$ is amenable. By [Sh99, Theorem 3.11], $P$ has a normal co-compact  subgroup $P_0$ which fixes every point in the support of $\xi$. Thus, if $a\in$ g$\setminus\{0\}$ is such that $p(a)\in $ P(g) is in the support of $\xi$, then there exists a homomorphism $\chi:P_0\rightarrow \Bbb R^*$ such that $\gamma a\gamma^{-1}=\chi(\gamma)a$, for all $\gamma\in P_0$. Now, $P_1=$ ker$(\chi)$ stabilizes $a$ and our assumption implies that $P_1$ is amenable. Since $\chi$ is continuous and $\Bbb R^*$ is amenable, we deduce that $P_0$ is amenable. Finally, as $P_0$ is co-compact in $P$, it follows that $P$ is amenable. 

Altogether, there exist an amenable subgroup $P<H$ and a   Borel function $\phi:X_0\rightarrow H/P$ such that $\phi(\theta(x))=w_{\theta}^1(x)\phi(x)$, for almost every $x\in X_0$, for all $\theta\in [R]$. 

Now, let us redefine $r:G\times G\rightarrow G$ as $r(x,y)=x^{-1}y$ and modify $\rho=p\circ q\circ r$ and $\pi$ accordingly. Repeating the above argument yields an amenable subgroup $Q<K$ and a  Borel function $\psi:X_0\rightarrow K/Q$ and  such that $\psi(\theta(x))=w_{\theta}^2(x)\psi(x)$, for almost every $x\in X_0$, for all $\theta\in [R]$ (note that $\tilde\eta$ is the weak limit of $p_1^*\mu_{n_k}$, so its support, $X_0$, does not depend on the definition of $\rho$).

Set $Y=H/P\times K/Q$ and $\tau:=(\phi,\psi):X_0\rightarrow Y$.
Then $\tau(\theta(x))=w_{\theta}(x)\tau(x)$, for almost every $x\in X_0$, for all $\theta\in [R]$.  Here on $Y=(H\times K)/(P\times Q)$ we consider the left multiplication action of $\Gamma\times\Lambda$.
 Since $P\times Q$ is amenable and $\Gamma\times\Lambda$ is discrete, this action  is topologically amenable, in the sense of [An02] (see the proof of [BO08, Theorem 5.4.1.]).

\vskip 0.05in

 Now, Proposition 3.6 in [Io10] gives that $R_{|X_0}$ is hyperfinite. For the reader's convenience we provide a self-contained argument proving that $R_{|X_0}$ is hyperfinite.  Fix a sequence $\{\theta_i\}_{i\geqslant 1}\subset [R]$ such that $R=\cup_{i\geq 1}\{(\theta_i(x),x|x\in X\}$. 
Define $w:R\rightarrow \Gamma\times\Lambda$ by $w(x,y)=w_{\theta_i}(y)$, where $i$ is the smallest integer with $x=\theta_i(y)$.
 Then $\tau:X_0\rightarrow Y$ satisfies $$\tau(x)=w(x,y)\tau(y),\hskip 0.07in\text{for}\hskip 0.03in \eta-\text{almost every}\hskip 0.05in (x,y)\in R_{|X_0}\tag d$$

 Since the action $\Gamma\times\Lambda\curvearrowright Y$ is topologically amenable,  by Connes-Feldman-Weiss' theorem ([CFW81]), its orbit equivalence relation $T$ is hyperfinite with respect to any measure on $Y$. Let $\eta_{|X_0}$ be the measure on $X_0$ given by $\eta_{|X_0}(A)=\eta(A\cap X_0)$. Then we can find an increasing sequence $T_n$ of finite equivalence relations on $Y$ such that $T=\cup_{n\geqslant 1}T_n$,  up to $\tau_*(\eta_{|X_0})$-null sets.

For every $n\geqslant 1$, set $R_n=\{(x,y)\in R_{|X_0}|(\tau(x),\tau(y))\in T_n\}$.  Then $R_n$ is an increasing sequence of subequivalence relations of $R_{|X_0}$ . By (d) we have that $\cup_{n\geqslant 1}R_n=R_{|X_0}$, up to $\eta$-null sets. Thus, to show that $R_{|X_0}$ is hyperfinite, it is enough to argue that $R_n$ is hyperfinite, for all $n\geqslant 1$. Now, if $R_0=\{(x,y)\in R_{|X_0}|\tau(x)=\tau(y)\}$, then $R_0$ has finite index in each $R_n$. Therefore, we can further reduce to proving that $R_0$ is hyperfinite.

Let us first prove this under the additional assumption that $R_0$ is ergodic. Then we can find $\tau\in Y$ such that $\tau(x)=\tau$, for almost every $x\in X_0$. Denote by $L$ the stabilizer of $\tau$ in $\Gamma\times\Lambda$ and by $S_0$ the orbit equivalence relation of the action $L\curvearrowright G$ (recall that $\Gamma\times\Lambda$ acts on $G$ by left-right multiplication). By using (d) we derive that $w(x,y)\in L$, for almost every $(x,y)\in R_0$. Thus, $R_0\subset {S_0}_{|X_0}$.
On the other hand, since $P,Q$ are amenable, it follows that $L$ is amenable. Thus, $S_0$ and $R_0$ are hyperfinite ([CFW81]).

If $R_0$ is not necessarily ergodic, we need to consider its ergodic decomposition. Let $\Cal M$ be the set of ergodic $R_0$-invariant probability measures on $X_0$ (viewed as a Borel subset as $\Cal M(X_0)$). Then there is an $R_0$-invariant Borel map $m:X_0\rightarrow \Cal M$ such that $\eta_{|X_0}=\int_{X_0}m(z) d\eta_{|X_0}(z)$ (see e.g. [KM04, Theorem 3.3]). Equation (d) implies that the set of $z\in X_0$ such that $\tau(x)=w(x,y)\tau(y),$ for $m(z)$-almost every $(x,y)\in R_0$, has full measure. Since $m(z)$ is $R_0$-ergodic, arguing as in the previous paragraph yields that $R_0$ is hyperfinite with respect to $m(z)$, for $\eta$-almost every $z\in X_0$. 
Since $\eta_{|X_0}=\int_{X_0}m(z) d\eta_{|X_0}(z)$, we conclude that $R_0$ is hyperfinite.

\vskip 0.05in

If $R_{|X_0}$ is {\it not ergodic}, the one proceeds as in the last paragraph. Consider the integral decomposition $\eta_{|X_0}=\int_{X_0}m(z) d\eta_{|X_0}(z)$, where $m(z)$ are ergodic $R_{|X_0}$-invariant probability measures. Then (c) (and the analogous identity for $w_{\theta}^2$ obtained by redefining $\rho$) holds when $\eta$ is replaced with $m(z)$, for almost every $z\in X_0$. Now, for such $z$, the above proof yields that $R_{|X_0}$ is $m(z)$-hyperfinite. Finally, this gives that $R_{|X_0}$ is $\eta_{|X_0}$-hyperfinite. 
\hfill$\square$ 
\vskip 0.1in

We are now ready to prove the results announced in the introduction.
\vskip 0.05in
\noindent
{\it Proof of Theorem A.} The first part is immediate by Theorem D. For the moreover part, recall Kazhdan's result: any connected semisimple Lie group $G$ with finite center whose simple factors have real-rank $\geqslant$ 2 has property (T) ([Ka67], see [Zi84, Theorem 7.4.2]).\hfill$\square$
\vskip 0.1in
\noindent
{\it Proof of Corollary B.} Denote by $\mu$ the probability measure on $X$. Let $(Y,\nu)$ be a pmp $\Gamma$-space and $p:X\times Y\rightarrow G/\Lambda$ be a measurable, quotient $\Gamma$-map.  If $A:=\{(x_1,x_2,y)\in X\times X\times Y|p(x_1,y)=p(x_2,y)\}$, then the conclusion is equivalent to $(\mu\times\mu\times\nu)(A)=1$. This implies that we may assume that the action $\Gamma\curvearrowright (Y,\nu)$ is ergodic.

Now, the action $\Gamma\curvearrowright (G/\Lambda,m_{G/\Lambda})$ is rigid by Theorem A. 
By [Io09, Proposition 3.3] we get that $A$  has positive measure. Since $\Gamma\cdot i$ is infinite for all $i\in I$, the action $\Gamma\curvearrowright (X,\mu)$ is weakly mixing. Hence, the product action of $\Gamma$ on $X\times X\times Y$ is ergodic. Since $A$ is invariant under this  action, the conclusion follows.
\hfill$\square$

\vskip 0.05in
{\it Remark.} In the case $Y$ is a single point space,  Corollary B has been first proved by Furman in  
[Fu07, Remark 1.15.(2)] by using entropy. When $\Gamma$ has property (T), this also follows from [Fu07, Theorem 1.14.] by using Popa's cocycle superrigidity theorem.
\vskip 0.1in

\noindent
{\it Proof of Theorem C.} Let $X\subset G$ be a fundamental set for the right $\Lambda$-action endowed with the probability measure $m(X)^{-1}m_{|X}$. Set $T=\{(x,y)\in X\times X|x\in \Gamma y\Lambda\}$.  Since the point stabilizers under the adjoint action of $G$ on its Lie algebra are amenable,  Theorem E implies that
 any subequivalence relation $R\subset T$, admits an $R$-invariant measurable partition $X=X_0\cup X_1$ such that $R_{|X_0}$ is hyperfinite and $R_{|X_1}$ is rigid.
 Since $\phi:G/\Lambda\ni x\Lambda\rightarrow x\Lambda\cap X\in X$ is a measure preserving isomorphism with $\phi(S)=T$, we are done.
 \hfill$\square$

\head References. \endhead
\item {[An02]} C. Anantharaman-Delaroche: {\it Amenability and exactness for dynamical systems and their C$^*$-algebras}, Trans. Amer. Math. Soc. {\bf 354} (2002),  4153--4178.
\item {[BO08]} N. P. Brown, N. Ozawa: {\it $C^*$-algebras and finite-dimensional approximations}, Graduate Studies in Mathematics, 88. American Mathematical Society, Providence, RI, 2008. xvi+509 pp.  
\item {[CJ85]} A. Connes, V. F. R. Jones: {\it Property (T) for von Neumann algebras,} Bull. London Math.
Soc. {\bf 17} (1985), 57--62.
\item {[CFW81]} A. Connes, J. Feldman, B. Weiss: {\it  An amenable equivalence relation is generated by a
single transformation}, Ergodic Theory Dynam. Systems {\bf 1} (1981), 431--450.
\item {[FM77]} J. Feldman, C. C. Moore: {\it Ergodic equivalence relations, cohomology, and von Neumann
algebras, II}, Trans. Amer. Math. Soc. {\bf 234} (1977), 325-–359.
\item {[Fu07]} A. Furman: {\it On Popa's Cocycle Superrigidity Theorem}, Inter. Math. Res. Notices IMRN,  {\bf 2007} (2007), 1--46, Art. ID rnm073.
\item {[GP05]} D. Gaboriau, S. Popa: {\it An uncountable family of non orbit equivalent actions of $\Bbb F_n$},
J. Amer. Math. Soc. {\bf 18} (2005), 547-–559.
\item {[Ga08]} D. Gaboriau: {\it Relative Property (T) Actions and Trivial Outer Automorphism Groups}, preprint arXiv:0804.0358.
\item {[Io07]} A. Ioana: {\it Orbit inequivalent actions for groups containing a copy of $\Bbb F_2$}, preprint  math/0701027.
\item {[Io09]} A. Ioana: {\it Non-orbit equivalent actions of $\Bbb F_n$}, Annales scientifiques de l'ENS {\bf 42}, fascicule 4 (2009), 675--696.
\item {[Io10]} A. Ioana: {\it Relative property (T) for the subequivalence relations induced by the action of
$\text{SL}(2,\Bbb Z)$ on $\Bbb T^2$}, Adv. Math. {\bf 224} (2010), 1589--1617.
\item {[Ka67]} D. Kazhdan: {\it On the connection of the dual space of a group with the structure of its closed subgroups}, Funct. Anal. and its Appl. {\bf 1}(1967), 63--65.
\item {[KM04]} A. S. Kechris, B. D. Miller: {\it Topics in orbit equivalence}, Lecture Notes in Mathematics, vol. 1852, Springer-Verlag, Berlin, 2004. 
\item {[Ma82]} G. Margulis: {\it Finitely-additive invariant measures on Euclidian spaces}, Ergodic Theory
Dynam. Systems 2(1982), 383-–396.
\item {[MvN36]} F. J. Murray, J. Von Neumann: {\it On rings of operators}, Ann. of Math. (2) {\bf 37} (1936),
no. 1, 116–-229.
\item {[Po06]} S. Popa: {\it On a class of type II$_1$ factors with Betti numbers invariants}, Ann. of Math. {\bf 163} (2006), 809--899.
\item {[PV10]} S. Popa and S. Vaes: {\it Actions of $\Bbb F_{\infty}$ whose II$_1$ factors and orbit equivalence relations have
prescribed fundamental group}, J. Amer. Math. Soc. {\bf 23} (2010), 383--403.
\item {[Sh99]} Y. Shalom: {\it Invariant measures for algebraic actions, Zariski dense subgroups and Kazhdan's property (T),} Trans. Amer. Math. Soc. {\bf 351}, 3387--3412.
\item {[Zi84]} R. J. Zimmer: {\it Ergodic theory and semisimple groups}, Monographs in Mathematics, 81. Birkhäuser Verlag, Basel, 1984. x+209 pp.  
\enddocument